\documentclass[12pt]{article}
\usepackage[russian]{babel}
\usepackage[cp1251]{inputenc}
\usepackage{amsmath,latexsym}
\usepackage[psamsfonts]{amssymb}
\usepackage[dvips]{graphicx}
\usepackage[mathcal]{euscript}
\begin{document}

УДК 517.968
\begin{center}
\textsc{ЭФФЕКТ ЕФИМОВА В МОДЕЛИ ФРИДРИХСА}
\end{center}
\begin{center}
\textbf{Ю. Х. Эшкабилов}
\end{center}

\begin{abstract}

Изучены дискретный спектр и эффект Ефимова в модели Фридрихса.
Приведены нобходимые и достаточные условия существования эффекта
Ефимова в модели Фридрихса. Изучен также дискретный спектр одного
модельного трехчастичного дискретного оператора Шредингера.
\end{abstract}

\textit{Ключевые слова:} модель Фридрихса, спектр, существенный
спектр, дискретный спектр, эффект Ефимова.\\

\textbf{1. Введение}\\

Пусть $u(x)$--вещественнозначная непрерывная функция на $[a,b]$,
$K$-компактный интегральный оператор в гильбертовом пространстве
$L_2[a,b]$.  Ряд вопросов квантовой механики и статистической
физики  (см. [1-3]) приводит к исследованию спектра оператора $H$
в гильбертовом пространстве $L_2[a,b],$ действующего по формуле
$$
H=U-K,\eqno (1)
$$
где $U$ -- оператор умножения на функцию $u(x)$, то есть
$Uf(x)=u(x)f(x),  f\in L_2[a,b]$.

Из классической теоремы Вейля о компактном возмущении, следует что
существенный спектр $\sigma_{e}(H)$ оператора $H$ состоит из
множества значений функции $u(x)$, т.е.
$\sigma_{e}(H)=[u_{\min},u_{\max}]$, где
$u_{\min}=\inf\limits_{x\in [a,\ b]} u(x),$
$u_{\max}=\sup\limits_{x\in [a,\ b]} u(x).$

Впервые оператор вида (1) был рассмотрен К. О. Фридрихсом [4] в
качестве простой модели теории возмущений существенного спектра.
Произвольный оператор вида (1) называется оператором в модели
Фридрихса.

 В работах Л. Д. Фаддеева [1] и О. А. Ладыженской, Л. Д.
Фаддеева [5] показано что, в случае  $u(x)=x$,  если ядро
интегрального оператора $K -$  гельдеровское ядро с показателем
$\mu
> \frac{1}{2},$ то оператор $U-K$ вне существенного спектра имеет
конечное число собственных значений конечной кратности.

Далее,  оператор в модели Фридрихса изучен в работах [6-9]. Пусть
$u(x)$-- вещественно-аналитическая фукнция в некоторой комплексной
окрестности отрезка $[a,b]$ и ядро оператора $K -$ симметричная и
аналитическая функция в некоторой окрестности квадрата $[a,b]^2$.
Доказано [6-9], что если число критических точек функции $u(x)$,
т.е. тех точек в которых $u^{'}(x)=0$ конечно и каждая из них
невырожденная, то дискретный спектр оператора $H$  конечен.

Вопрос о бесконечности числа собственных значений в модели (1),
лежащих вне существенного спектра является мало изученным , хотя с
математической точки зрения в этой модели можно ожидать  явление
эффекта Ефимова.

В настоящей заметке рассматривается вопрос о существовании
бесконечного числа собственных значений  в  модели Фридрихса.
Найдены необходимые и достаточные условия бесконечности
дискретного спектра в модели Фридрихса и как следствие полученного
результата доказана теорема о дискретном спектре "трехчастичного"
дискретного оператора. \ \ \ \ \ \ \ \ \ \ \ \\ \ \ \ \ \

\textbf{2. Некоторые обозначения и необходимые сведения}\\

Пусть $\mathcal{H}$ -- гильбертово пространство,
$A:\mathcal{H}\rightarrow\mathcal{H}$ -- линейный ограниченный
самосопряженный оператор.  Через $\sigma(A),$ $\sigma_e(A)$ и
$\sigma_d(A)$ обозначим, соответственно, спектр, существенный
спектр и дискретный спектр оператора $A$ [10]. Определим
вещественное число
$$
S_{\min}=S_{\min}(A)=\inf\limits_{\|x\|=1}(Ax,x) \mbox{ и }
S_{\max}=S_{\max}(A)=\sup\limits_{\|x\|=1}(Ax,x)
$$
Тогда $\sigma(A)\subset[S_{\min}, S_{\max}],$ при этом $S_{\min},\
S_{\max}\in \sigma(A)$.

Введем также следующие обозначения
$$
\mathcal{E}_{\min}=\mathcal{E}_{\min}(A)=\inf\{\lambda:\
\lambda\in \sigma_e(A)\},
$$
$$
\mathcal{E}_{\max}=\mathcal{E}_{\max}(A)=\sup\{\lambda:\
\lambda\in \sigma_e(A)\}.
$$

Число $\mathcal{E}_{\min}\in \sigma(A)$ ($\mathcal{E}_{\max}\in
\sigma(A)$) будем называть нижним краем (верхним краем)
существенного спектра оператора $A$. Например, для компактного
самосопряженного оператора $K:\mathcal{H}\rightarrow\mathcal{H}$
имеем $\sigma_{e}(K)=\{0\}$ и отсюда
$\mathcal{E}_{\min}(K)=\mathcal{E}_{\max}(K)=0$.

\textbf{Предложение 2.1.} \textit{(а) если
$S_{\min}=\mathcal{E}_{\min}$ ($S_{\max}=\mathcal{E}_{\max}$), то
у оператора $A$ отсутствует точка дискретного спектра, лежащяя
ниже (выше) нижнего (верхнего) края существенного спектра, т. е.
$\sigma_d(A)\cap(-\infty,\ \mathcal{E}_{\min}]=\emptyset$
($\sigma_d(A)\cap[\mathcal{E}_{\max},\ +\infty)=\emptyset$);}

(б) \textit{если $S_{\min}<\mathcal{E}_{\min}$
($S_{\max}>\mathcal{E}_{\max}$), то для оператора $A$ множество
точек дискретного спектра, лежащих ниже (выше) нижнего (верхнего)
края существенного спектра является непустым, т. е.
$\sigma_d(A)\cap(-\infty,\ \mathcal{E}_{\min}]\not=\emptyset$
($\sigma_d(A)\cap[\mathcal{E}_{\max},\ +\infty)\not=\emptyset$).}

Для заданного ограниченного самосопряженного оператора
$A:\mathcal{H}\rightarrow\mathcal{H}$ построим ограниченную
возрастающую последовательность $\mu_n=\mu_n(A),$ $n\in
\mathbb{N},$ следующим способом:
$$
\mu_1(A)=S_{\min}(A)=\inf\limits_{x\in \mathcal{H},\ \|x\|=1}(Ax,\
x).
$$
Допустим, что $\mu_1(A)=(Ax_1, x_1)$ для некоторого $x_1\in
\mathcal{H}, \|x_1\|=1$.

Определим
$$
\widetilde{\mu}_2(A)=\mu_1(A_1')=S_{\min}(A_1')=\inf\limits_{x\in
^\bot\mathcal{H}_1,\ \|x\|=1}(A_1'x, x),
$$
где $^\bot\mathcal{H}_1=\{x\in \mathcal{H}: (x, x_1)=0\}$ и $A_1'$
-- сужение оператора $A$ на подпространство $^\bot\mathcal{H}_1$.

Аналогично определим
$$
\widetilde{\mu}_k(A)=\mu_1(A_{k-1}')=S_{\min}(A_{k-1}')=\inf\limits_{x\in
^\bot\mathcal{H}_{k-1},\ \|x\|=1}(A_{k-1}'x, x),
$$
где $S_{\min}(A_{k}')=(A_k'x_{k+1},\ x_{k+1}),$ $\|x_{k+1}\|=1,$
$^\bot\mathcal{H}_k=\{x\in \mathcal{H}: (x,\ x_1)=...=(x,\
x_k)=0\},$ $A_k'$ -- сужение оператора $A$ на подпространство
$^\bot\mathcal{H}_k,$ $k\in \mathbb{N}.$

Теперь определим последовательность
$\{\mu_n(A)\}\subset\mathbb{R}$ по правилу
$$
\mu_n(A)=\sup\{\mu_1(A),\widetilde{\mu}_2(A),...,\widetilde{\mu}_n(A)\},
\qquad n\in \mathbb{N}. \eqno (2)
$$

Для заданного ограниченного самосопряженного оператора
$A:\mathcal{H}\rightarrow\mathcal{H}$ построим ограниченную
убывающую последовательность $\eta_n=\eta_n(A),$ $n\in N$
следующим способом:
$$
\eta_1(A)=S_{\max}(A)=\sup\limits_{y\in \mathcal{H}, \|y\|=1}(Ay,
y).
$$

Допустим, что $\eta_1(A)=(Ay_1, y_1)$ для некоторого $y_1\in
\mathcal{H}$, $\|y_1\|=1$. Определим число
$$\widetilde{\eta}_{2}(A)=\eta_1(A''_1)=S_{\max}(A''_1)=\sup\limits_{y\in \mathcal{H}_{1}^{\bot},
\|y\|=1}(A_{1}^{''} y, y).$$ где $\mathcal{H}^\bot_1=\{y\in
\mathcal{H}:\ (y, y_1)=0\}$ и $A_1''$ -- сужение оператора $A$ на
подпространство $\mathcal{H}^\bot_1.$

Аналогично определим
$$
\widetilde{\eta}_k(A)=\eta_1(A_{k-1}'')=S_{\max}(A_{k-1}'')=\sup\limits_{y\in
\mathcal{H}^\bot_{k-1}, \|y\|=1}(A_{k-1}''y, y),
$$
где $S_{\max}(A_k'')=(A_{k+1}''y_{k+1},\ y_{k+1}),$
$\|y_{k+1}\|=1,$ $\mathcal{H}^\bot_{k}=\{y\in \mathcal{H}:\ (y,
y_1)=...=(y, y_k)=0\},$ $A_{k}''$ -- сужение оператора $A$ на
подпространство $\mathcal{H}^\bot_{k},$ $k\in \mathbb{N}.$

Определим последовательность $\{\eta_n(A)\}\subset\mathbb{R}$ по
правилу
$$
\eta_n(A)=\inf\{\eta_1(A),\ \widetilde{\eta}_2(A),...,\
\widetilde{\eta}_n(A)\},\qquad n\in \mathbb{N}. \eqno(3)
$$

В книге [11] пользуясь методом принципа минимакса изучен
дискретный спектр  заданного  самосопряженного и ограниченного
снизу оператора  действующего в гильбертовом пространстве. В
частности, доказано существование собственных значений лежащих
ниже существенного спектра и бесконечность дискретного спектра для
некоторых многочастичных гамильтонианов.

Надо подчеркнуть, что для самосопряженных и ограниченных снизу
операторов существенный спектр $\sigma_{e}(\cdot)$ по построению
имеет удобный вид, т.е. он представляется только одним бесконечным
отрезком $[a,\infty)$. Для таких операторов дискретний спектр
лежит только в одном полуинтервале, т.е. существует число $c<a$
такое, что $\sigma_{d}(\cdot)\subset [c,a)$. Поэтому принцип
минимакса для таких операторов успешно применяется [11].

Существенный спектр ограниченного самосопряженного оператора может
состоять из нескольких отрезков, между этими отрезками могут
существовать точки дискретного спектра. Например, если для
заданного ограниченного самосопряженного оператора $ A $ \ \
$\sigma_{e}(A)=\bigcup\limits_{k}[a_k,b_k]$, где $b_k<a_{k+1},$ то
$\sigma_{d}(A)\subset[\mu_{1}(A),\mathcal{E}_{\min})\cup(\bigcup\limits_{k}
(b_k,a_{k+1}))\cup(\mathcal{E}_{\max},\eta_{1}(A)]$ , где
$\mathcal{E}_{\min}=\inf \{a_{1},a_{2},...\},
\mathcal{E}_{\max}=\sup \{b_{1},b_{2},...\}$. Принцип минимакса и
максимина нам позволяет изучать дискретные спектры самосопряженных
ограниченных операторов лежащих ниже нижнего края существенного
спектра $\mathcal{E}_{\min}$ и лежащих выше верхнего края
существенного спектра $\mathcal{E}_{\max}$.

 \textbf{Теорема 2.2} (принцип минимакса для
ограниченных операторов в операторной форме). \textit{Пусть
$A:\mathcal{H}\rightarrow\mathcal{H}$ -- ограниченный
самосопряженный
оператор. Тогда для каждого фиксированного $n\in \mathbb{N}$:\\
либо\\
а) существует $n$ собственных значений (считая вырожденные
собственные значения столько раз, какова их кратность), лежащие
ниже нижнего края $\mathcal{E}_{\min}(A)$ существенного спектра, а
$\mu_n(A)$ (2) есть $n$-е собственное
значение оператора $A$ (с учетом кратности),\\
либо \\
б) $\mu_n(A)$ -- нижний край существенного спектра, т. е.
$\mu_n(A)=\mathcal{E}_{\min}(A)$, при этом
$\mu_n=\mu_{n+1}=\mu_{n+2}=...$ и существует самое большее $n-1$
собственных значений (с учетом кратности), лежащих  ниже нижнего
края $\mathcal{E}_{\min}(A)$ существенного спектра оператора $A.$}

\textsl{Доказательство} проводится аналогично доказательству
теоремы о принципе минимакса для операторов ограниченных снизу
[11].

Теорема о принципе максимина для ограниченных операторов
формулируется следующим образом.

\textbf{Теорема 2.3} (принцип максимина для ограниченных
операторов в операторной форме ). \textit{Пусть
$A:\mathcal{H}\rightarrow\mathcal{H}$ -- ограниченный
самосопряженный оператор. Тогда для каждого
фиксированного $n\in \mathbb{N}$:\\
либо \\
а) существует $n$ собственных значений (считая вырожденные
собственные значения столько раз, какова их кратность), лежащие
выше верхнего края $\mathcal{E}_{\max}(A)$ существенного спектра,
а $\eta_n(A)$ (3) есть $n$-е собственное
значение оператора $A$ (с учетом кратности),\\
либо \\
 б) $\eta_n(A)$ -- верхний край существенного спектра, т. е.
$\eta_n(A)=\mathcal{E}_{\max}(A),$ при этом
$\eta_n=\eta_{n+1}=\eta_{n+2}=...$ и существует самое большее
$n-1$ собственных значений (с учетом кратности), лежащих выше
верхнего края $\mathcal{E}_{\max}(A)$ существенного спектра
оператора $A.$}

Линейный ограниченный оператор $A: \mathcal{H} \rightarrow
\mathcal{H}$ называется положительным и пишется $A\geq 0$ или
$0\leq A,$ если $(Ax, x) \geq 0, \forall x \in \mathcal{H}.$

\textbf{Следствие 2.4.} \emph{Пусть $A, B: \mathcal{H} \rightarrow
\mathcal{H}$ -линейные ограниченные самосопряженные операторы и $A
\leq B.$ Тогда
 $$\mbox{(а) Если}\
\mathcal{E}_{min}(A)=\mathcal{E}_{min}(B),\  \mbox{то} \
\mu_n(A)\leq \mu_n(B), n\in \mathbb{N} \eqno(4)$$
$$ \mbox{(б) Если}\ \mathcal{E}_{max}(A)=\mathcal{E}_{max}(B), \mbox{то} \ \eta_n(A)\leq
\eta_n(B), n\in \mathbb{N} \eqno(5) $$}

Будем говорить, что точка $\lambda\in \sigma_e(A)$ изолирована
слева (справа) на существенном спектре оператора $A$, если
существует число $\delta>0$ такое, что $(\lambda-\delta,\
\lambda]\cap\sigma_e(A)=\{\lambda\}$ ($[\lambda,\
\lambda+\delta)\cap\sigma_e(A)=\{\lambda\}$).

\textbf{Определение 2.1.} Пусть точка $\lambda\in \sigma_e(A)$
изолирована слева на существенном спектре оператора $A$. Если
существует последовательность $\{\lambda_n\}$ собственных значений
оператора $A$, удовлетворяющая условиям $\lambda_n<\lambda$ для
всех $n\in \mathbb{N}$ и
$\lim\limits_{n\rightarrow\infty}\lambda_n=\lambda$, то говорят,
что в точке $\lambda$ \emph{слева существует эффект Ефимова.}

\textbf{Определение 2.2.} Пусть точка $\lambda\in \sigma_e(A)$
изолирована справа на существенном спектре оператора $A$. Если
существует последовательность $\{\lambda_n\}$ собственных значений
оператора $A$, удовлетворяющая условиям $\lambda_n>\lambda$ для
всех $n\in \mathbb{N}$ и
$\lim\limits_{n\rightarrow\infty}\lambda_n=\lambda$, то говорят,
что в точке $\lambda$ \emph{справа существует эффект Ефимова.}\\

{\bf 3. Модель Фридрихса}\\

Пусть $u(t)-$ неотрицательная непрерывная функция на $[0,1]$ и
$u^{-1}(\{0\})\bigcap [0,1]\neq \emptyset; \,\,\ k(t,s)\in
L_2[0,1]^2$ и $k(t,s)=\overline{k(s,t)}.$ В пространстве
$L_2[0,1]^2$ рассмотрим самосопряженный оператор $$H=U-K,
\eqno(6)$$ где
$$(U f)(t)=u(t)f(t), \ \ f\in L_2[0,1]$$
$$(Kf)(t)=\int_0^1k(t,s) f(s)d \mu(s), \ \ f \in L_2[0,1].$$

Из теоремы Вейля [11] о компактном возмушении существенного
спектра получаем, что $\sigma_e(H)=\sigma(U)=[0,u_{max}]$ и отсюда
$\mathcal{E}_{min}(H)=0, \mathcal{E}_{max}(H)=u_{max}.$

\textbf{Теорема 3.1.} \textit{Если интегральный оператор $K$ в
модели (6) имеет конечное число положительных(отрицательных)
собственных значений, то оператор $U-K$ также имеет конечное число
собственных значений лежащих ниже(выше) существенного спектра.}

\textit{Доказательство} а) Пусть оператор $K$ имеет только $n-1$
число положительных собственных значений. Тогда оператор $-K$
также имеет ровно $n-1$ отрицательных собственных значений. Имеем
$U\geq 0.$ Тогда из $\mathcal{E}_{min}(-K)=\mathcal{E}_{min}(U-K)$
в силу неравенства (4) получаем, что $\mu_k(-K) \leq \mu_k(U-K),
\forall k \in \{1,2,\dots,n\}.$ В силу принципа минимакса
$\mu_n(-K)=\mu_{n+1}(-K)=\dots=0.$ Следовательно, оператор $U-K$
ниже существенного спектра имеет не больше $n-1$  собственных
значений.

б) Пусть оператор $K$ имеет только $n-1$  отрицательных
собственных значений. Тогда оператор $u_{max}\cdot E-K$ также
имеет только $n-1$ собственных значений лежащих выше числа
$u_{max}$, где $E-$ тождественный оператор в $L_2[0,1]$. Имеем
$\mathcal{E}_{max}(U-K) = \mathcal{E}_{max}(u_{max}E-K).$ Тогда из
соотношений $u_{max} E-U\geq 0$ получаем, что
$$\eta_k(U-K)\leq \eta_k(u_{max}E-K), \ \forall k \in \{1,2,\dots, n\}.$$

Однако в силу принципа максимина имеем
$u_{max}=\eta_n(u_{max}E-K)=\eta_{n+1}(u_{max}E-K)=\dots\ .$
Отсюда $\eta_{n+k}(U-K)=u_{max}, k=0,1,\dots,$ т.е. оператор $U-K$
имеет не более $n-1$ собственных значений, лежащих выше
существенного спектра. $\blacksquare$

\textbf{Следствие 3.2.} \textit{Чтобы у оператора $U-K$ ниже
(выше) существенного спектра существовал эффект Ефимова
необходимо, чтобы компактный оператор $K$ имел бесконечное число
положительных (отрицательных) собственных значений.}

\textbf{Теорема 3.3.} \textit{(а) Чтобы у оператора $U-K$ ниже
существенного спектра существовал эффект Ефимова необходимо, чтобы
$$\lim_{n\rightarrow \infty} \int_0^1 u(t)|\varphi_n(t)|^2
d\mu(t)=0,$$ где $\{\varphi_n\}_{n \in \mathbb{N}}-$ система
ортонормированных собственных функций соответствующих собственным
значениям оператора $U-K$ лежащих ниже существенного спектра;
\\(б) Чтобы у оператора $U-K$ выше сушественного спектра
существовал эффект Ефимова необходимо, чтобы
$$\lim_{n\rightarrow \infty} \int_0^1 u(t)|\psi_n(t)|^2
d\mu(t)=u_{max},$$ где $\{\psi_n\}_{n \in \mathcal{N}}-$ система
ортонормированных собственных функций соответствующих собственным
значениям оператора $U-K$ лежащих выше сушественного спектра.}

\textit{Доказательство.} а) Пусть для оператора $U-K$ в точке
$\lambda =0$ слева существует эффект Ефимова. Тогда существует
последовательность отрицательных собственных значений оператора
$U-K,$ т.е. $\mu_n(U-K)<0, \forall n \in \mathbb{N}.$  Пусть
$\{\varphi_n\}_{n\in \mathbb{N}}$ - система ортонормированных
собственных функций соответствующих собственным значениям
$\mu_n(U-K), n \in\mathbb{N}.$ Имеем $(U\varphi_n,
\varphi_n)-(K\varphi_n, \varphi_n) = \mu_n(U-K) < 0, n \in
\mathbb{N}.$ Отсюда получаем, что $0 \leq (U\varphi_n, \varphi_n)
< (K\varphi_n, \varphi_n), n \in \mathbb{N}.$ В силу компактности
оператора $K$ имеем $\lim_{n \rightarrow \infty} (K\varphi_n,
\varphi_n)=0.$ Следовательно, $\lim_{n  \rightarrow \infty}
(U\varphi_n, \varphi_n)=0,$ т.е.
$$\lim_{n\rightarrow \infty} \int_0^1 u(t)|\varphi_n(t)|^2
d\mu(t)=0.$$

б) Пусть для оператора $U-K$ в точке $\lambda =u_{max}$ справа
существует эффект Ефимова. Тогда существует последовательность
$\lambda_n$ собственных значений оператора $U-K,$ удовлетворяющая
неравенству $\lambda_n > u_{max}, \forall n \in \mathbb{N},$ т.е.
$\eta_n(U-K)>u_{max}. $ Пусть $\{\psi_n\}_{n\in \mathbb{N}}$
-система ортонормированных собственных функций соотвествующих
собственным значениям $\eta_n(U-K), n \in \mathbb{N}.$ Тогда имеем
$(U\psi_n, \psi_n)-(K\psi_n, \psi_n) = \eta_n(U-K) > u_{max}, n
\in \mathbb{N}.$  Отсюда  $u_{max} + (K\psi_n, \psi_n) < (U\psi_n,
\psi_n) \leq u_{max}, n \in \mathbb{N}.$ Следовательно,
$$\lim_{n\rightarrow \infty} (U\psi_n, \psi_n) = \int_0^1 u(t)|\psi_n(t)|^2
d\mu(t)=u_{max}. \eqno \blacksquare$$

\textbf{Теорема 3.4.} \textit{Если существует ортонормированная
система $\{f_n\}_{n \in \mathbb{N}}\subset L_2 [0,1]$
удовлетворяющая условию
$$ \int_0^1 u(t)|f_k(t)|^2
d\mu(t) <  c_k, \forall k \in \mathbb{N}, \eqno (7) $$ где
$c_k=(Kf_k,f_k),$ тогда слева от точки $\lambda =0$ спектра
оператора $U-K$ существует эффект Ефимова.}

\textit{Доказательство.} Пусть в пространстве $L_2[0,1]$
существует ортонормированная система $\{f_n\}, \ n \in \mathbb{N}$
удовлетворяющая условию (7). Предположим, что система $\{f_n\}, n
\in \mathbb{N}$ упорядочена в следующим смысле: $(Hf_k, f_k) \leq
(Hf_{k+1}, f_{k+1}), k \in \mathbb{N}.$ Имеем
$$S_{min}(H)=\mu_1(H) \leq (H f_1, f_1) < 0,$$
т.е. существует $\varphi_1 \in L_2 [0,1], \|\varphi_1\|=1$ такая,
что $\mu_1(H)=(H\varphi_1, \varphi_1)<0$ и по теореме о принципе
минимакса $\mu_1(H)$ является собственным значением оператора
$H=U-K.$ Для каждого $k \in \mathbb{N}$ определим подпростанство
$^{\perp}L_2^{(k)}[0,1]:$
$$^\perp L_2^{(k)}[0,1]=\{f \in L_2[0,1]: (f, f_j)=0,\,\ j=1,2,\dots, k\}.$$

Пусть $H_k$ сужение оператора $H$ на подпространство $^\perp
L_2^{(k)}[0,1], k \in \mathbb{N}.$ Для оператора $H_1$ имеем
$$\mu_1(H_1)=S_{min}(H_1)=
\inf \limits_{f \in ^\perp L_2^{(1)}[0,1], \|f\|=1}(H_1f, f)
\leq$$
$$\leq(H_1f_2, f_2)=(Hf_2, f_2)<0$$

Так как $f_2 \in  {^\perp} L_2^{(1)}[0,1], \|f_2\|=1,$ то
существует $\varphi_2 \in  {^\perp} L_2^{(1)}[0,1],
\|\varphi_2\|=1$ такая, что $\mu_1(H_1)=(H_1\varphi_2,
\varphi_2)=(H\varphi_2, \varphi_2)<0.$ Для каждого  $k \geq 2$
повторяя аналогичное рассуждение имеем
$$\mu_1(H_k)=S_{min}(H_k)=\inf \limits_{f \in {^\perp} L_2^{(k)}[0,1], \|f\|=1}
(H_kf, f)\leq $$ $$\leq (H_kf_{k+1}, f_{k+1})=(Hf_{k+1},
f_{k+1}).$$

Так как $f_{k+1} \in {^\perp} L_2^{(k)}[0,1], \|f_{k+1}\|=1, \ k
\geq 2,$ то существует $\varphi_{k+1} \in$ $
^{\perp}L_2^{(k)}[0,1], \|\varphi_{k+1}\|=1$ такая, что
$\mu_1(H_k)=(H_k\varphi_{k+1}, \varphi_{k+1})= (H\varphi_{k+1},
\varphi_{k+1})<0,$ где $k \geq 2.$ Таким образом, для каждого $k
\in \mathbb{N}$ оператор сужения $H_k$ имеет собственное значение
$\mu_1(H_k)<0$ и $\mu_1(H_k)\leq \mu_1(H_{k+1}), k\in \mathbb{N}.$
Следовательно, каждое число $\omega_k = \mu_1(H_k)<0, k=0, 1, 2,
\dots$ является собственным значением оператора $H,$ где $H_0=H$ и
$\lim \limits_{n\rightarrow \infty} \omega_n=0.$ $\blacksquare$

\textbf{Теорема 3.5.} \textit{Если сушествует ортонормированная
система $\{f_n\} \subset L_2[0,1]$ удовлетворяющая условию
$$\int_0^1u(t)|f_k(t)|^2d\mu(t)> c_k+ u_{max},\,\ \forall k \in \mathbb{N} \eqno (8)$$
где $c_k=(Kf_k,f_k)$, тогда справа от точки $\lambda=u_{max}$
спектра оператора $U-K$ существует эффект Ефимова.}

\textbf{Замечание 3.6.} \textit{Если в модели Фридрихса (6)
функция $u(t)$ не удовлетворяет условию
$u^{-1}(\{0\})\bigcap[0,1]\neq \emptyset,$ то сдвигая оператор $H$
на $-u_{min}$ имеем оператор $\widetilde{H}=\widetilde{U}-K$ в
модели Фридрихса, где $(\widetilde{U}f)(t)=(u(t)-u_{min})f(t), f
\in L_2[0,1]$ и функция $h(t)=u(t)-u_{min}$ удовлетворяет условию
$h(t)\geq 0, \ t \in [0,1]$ и $h^{-1}(\{0\}) \bigcap [0,1]\neq
\emptyset.$ Следовательно, получая информацию о спектре оператора
$\widetilde{H}$ мы получаем  полную информацию о спектре оператора
$H,$ так как: если  $\widetilde{\omega}$ является $m-$ кратным
собственным значением оператора $\widetilde{H},$ то число
$\widetilde{\omega}+u_{min}$ является $m-$ кратным собственным
значением оператора $H$ и наоборот, при этом
$\sigma_e(H)=[u_{min}, u_{max}].$}
\\

\textbf{4. "Трехчастичный" \ оператор Шредингера $T$}\\

Одним из замечательных результатов в спектральной теории
многочастичных гамильтонианов является  эффект Ефимова. Этот
эффект для трехчастичных гамильтонианов впервые был обнаружен в
1970 году В.Н.~Ефимовым [12]. С тех пор эта задача изучалась во
многих работах (см., например, [13-18]). В работах [19-22] изучен
эффект Ефимова для трехчастичного дискретного оператора
Шредингера. В этом параграфе мы рассмотрим модель трехчастичного
дискретного оператора $T$ и изучим условия существования эффекта
Ефимова.

Пусть $T$ ограниченный самосопряженный оператор в  $L_2[0,1]^2$
 заданный по правилу
$$T=T_0-(T_1+T_2), \eqno (9)$$
где $$(T_0f)(x,y)=(u(x)+v(y))f(x,y), \  f \in L_2[0,1]^2$$
$$(T_1f)(x,y)=\int\limits_0^1 \int\limits_0^1k_1(x,s)\delta(y-t) f(s,t)
d\mu(s)d\mu(t), \ f \in L_2[0,1]^2$$
$$(T_2f)(x,y)=\int\limits_0^1 \int\limits_0^1k_2(y,t)\delta(x-s) f(s,t)
d\mu(s)d\mu(t), \ f \in L_2[0,1]^2$$ Здесь $u(x), v(x)-$ заданные
неотрицательные непрерывные функции на $[0,1]$ и
$u^{-1}(\{0\})\bigcap [0,1]\neq \emptyset, v^{-1}(\{0\})\bigcap
[0,1]\neq \emptyset;$ $k_1(x,s), k_2(x,s)-$ заданные непрерывные
функции на квадрате $[0,1]^2$ и $k_1(x,s)=\overline{k_1(s,x)},$
$k_2(x,s)=\overline{k_2(s,x)},  \delta(x)-$ дельта функция Дирака
на $[0,1].$

Операторы $T_1$ и $T_2$ не являются компактными. Существенный
спектр оператора в модели (9) в более общем случае изучен в работе
автора [23]. В модели Хаббарда [23] возникает "трехчастичный" \
дискретный оператор, который представляется в виде (9). Оператор
$T$ (9) является общим аналогом трехчастичного дискретного
оператора, соответствующего системе состоящей из двух свободных
электронов и одной примеси на решетке, и его в дальнейшим будем
называть "трехчастичным" \ дискретным оператором Шредингера.

\textbf{Замечание 4.1.} \textit{Оператор $T$ (9) можно
представить, как оператора в модели Фридрихса с некомпактым ядром,
т.е. $T=T_0 -W,$ где $W$ интегральный оператор в $L_2[0,1]^2$ с
некомпактным ядром $q(x,y;
s,t)=k_1(x,s)\delta(y-t)+k_2(y,t)\delta(x-s).$}

Наряду с гамильтонианом $T$ мы будем рассматривать два оператора
$H_1$ и $H_2$ в пространстве $L_2[0,1],$  каждый из которых
представляется в модели Фридрихса в следующим виде:
$$H_1 = U_1 -K_1, \eqno(10)$$
где

$$(U_1{\varphi})(x)=u(x)\varphi(x), \ \varphi \in L_2[0,1],$$
$$(K_1{\varphi})(x) =\int\limits_0^1 k_1(x,s) \varphi(s)d\mu(s), \
\varphi \in L_2[0,1]$$ и
$$H_2 = U_2 -K_2, \eqno(11)$$
где $$(U_2{\varphi})(x)=v(x)\varphi(x), \ \varphi \in L_2[0,1],$$
$$(K_2{\varphi})(x)=\int\limits_0^1 k_2(x,s) \varphi(s)d\mu(s), \ \varphi \in L_2[0,1].$$

В дальнейшем мы предположим, что ядра $k_1(x,s)$ и $k_2(x,s)$ в
модели (9) удовлетворяют условиям $K_1 \geq 0$ и $K_2 \geq 0,$
т.е. $K_1,\, K_2 -$ положительные операторы.

Обозначим через $\alpha_1, \alpha_2, \dots, \alpha_k, \dots$
отрицательные собственные значения оператора $H_1$, через $\rho_1$
обозначим их количество и предположим, что $\alpha_k \leq
\alpha_{k+1}$. Тогда имеем
$\sigma_d(H_1)=\{\alpha_k\}_{k=1}^{\rho_1}.$ Аналогично, через
$\beta_1, \beta_2,\dots, \beta_k, \dots$ обозначим отрицательное
собственное значение оператора $H_2,$ через $\rho_2$ их количество
и предположим, что $\beta_k \leq \beta_{k+1}$. Тогда
$\sigma_d(H_2)=\{\beta_k\}_{k=1}^{\rho_2}.$

Пользуясь тензорным произведением пространств и операторов [10],
гамильтониан $T$ мы можем представить в виде $T=H_1\otimes E +
E\otimes H_2.$   Для спектра оператора $T$ имеем [10]:
$\sigma(T)=\sigma(T_1)+\sigma(T_2).$ Отсюда следует [23]:

\textbf{Теорема 4.1.} а) $\sigma_e(T)=\sigma_0\cup\sigma_1\cup
\sigma_2,$ где $\sigma_0=[0,u_{max}+v_{max}], u_{max}=\sup u(x),
v_{max} = \sup v(x)$ и
$$\sigma_1 = \overline{\bigcup\limits_{k=1}^{\rho_1}[\alpha_k,
v_{max} + \alpha_k]},\ \ \sigma_2 = \overline{\bigcup
\limits_{k=1}^{\rho_2}[\beta_k, u_{max} + \beta_k]}.$$

б) $\sigma_d(T)=\{\omega: \omega= \alpha+\beta \notin
\sigma_1\cup\sigma_2, \alpha \in \sigma_d (H_1), \beta \in
\sigma_d (H_2)\}.$

\ \\ \ \ \ \ Легко проверить, что $\sigma(T) \subset[\alpha_1 +
\beta_1, u_{max}+v_{max}],$ при этом $\alpha_1+ \beta_1 \in
\sigma_d(T).$

\textbf{Следствие 4.2.} \textit{Существенный спектр гамильтониана
$T$ состоит из объединения конечного число отрезков.}

\textit{Доказательство.} Если дискретный спектр операторов $H_1$ и
$H_2$ конечен  ({\bf случай 1}), то из утверждения (а) теоремы
4.1. следует, что $\sigma_e(T)$ состоит из объединения конечного
числа отрезков.

\textbf{Случай 2:} Допустим, что $\sigma_d(H_1)-$ бесконечен и
$\sigma_d(H_2)-$ конечен. Тогда для возрастающей
последовательности $\alpha_n, n \in \mathbb{N}$ имеем
$\lim\limits_{n \rightarrow \infty}\alpha_n=0.$ Следовательно,
существует натуральное число $n_0=n_0(H_1)$ удовлетворяшее
условию: $\alpha_{n+1}-\alpha_n < v_{max}, \forall n \geq n_0.$

Для натурального числа $n_0=n_0(H_1)$ имеем
$$\bigcup\limits_{k=n_0}^{\infty}[\alpha_k, v_{max}+\alpha_k]=
[\alpha_{n_0}, v_{max}].$$ Следовательно,
$$\sigma_1 = \bigcup\limits_{k=1}^{n_0-1}[\alpha_k,
v_{max}+\alpha_k]\bigcup [\alpha_{n_0}, v_{max}],$$ т.е.
$\sigma_1$ состоит из объединения конечного числа отрезков, а
множество $\sigma_2$ по предположению обладает этим свойствам.
Значит $\sigma_e(T)$ представляется в виде объединения конечного
число отрезков.

\textbf{Случай 3:} $\sigma_d(H_1)-$ конечен и $\sigma_d(H_2)-$
бесконечен. Тогда опять, как в случае 2 сушествует натуральное
число $m_0 = m_0(H_2)$ такое, что $\beta_{n+1}-\beta_n < u_{max},
\forall n \geq m_0.$

 Следовательно, имеем
$$\sigma_2 = \bigcup\limits_{k=1}^{m_0-1}[\beta_k,
u_{max}+\beta_k] \bigcup [\beta_{m_0}, u_{max}].$$

\textbf{Случай 4:} $\sigma_d(H_1)$ и $\sigma_d(H_2)-$ бесконечены.
Тогда определяя натуральные числа $n_0=n_0(H_1)$ и $m_0=m_0(H_2),$
как в предыдуших случаях, мы будем множество $\sigma_e(T)$
предcтавлять в виде объединения конечного числа отрезков.
$\blacksquare$

\textbf{Следствие 4.3.} \textit{Всякая изолированная точка в
спектре оператора $T$ является конечнократным собственным
значением.}

\textbf{Следствие 4.4.} \textit{\\(а) если $\sigma_d
(H_1)=\emptyset$ или $\sigma_d(H_2)=\emptyset,$ то $\sigma_d(T)=
\emptyset;$\\
б) если $|\sigma_d(H_1)| > 0$ и $|\sigma_d(H_2)| > 0,$ то
$|\sigma_d(T)| > 0,$ где  $|\Omega|-$ мощность множества $\Omega
.$\\ (в) если $|\sigma_d(H_1)| =p$ и $|\sigma_d(H_2)| = q,$ то
$|\sigma_d(T)| \leq pq.$}

Обозначим через $n_e(T)$ количество взаимно непересекающихся
отрезков, объединением которых является множество $\sigma_e(T).$
Таким образом существенный спектр $\sigma_e(T)$ представляется в
виде

$$\sigma_e(T)=\bigcup\limits_{k=1}^{n_e(T)}[S_k,S_k'], \ S_k <
S_k' < S_{k+1}. \eqno (12)$$

Нетрудно проверить, что в (12) $S_1= min \{\alpha_1, \beta_1 \},$
при $|\sigma_d(H_1)| > 0$ и $|\sigma_d(H_2)| > 0.$

\textbf{Теорема 4.5.} \textit{Пусть 0<$|\sigma_d(H_1)|=p$ и  $
|\sigma_d(H_2)|=\infty$. Если для $\alpha \in \sigma_d(H_1)$
найдется число $S_k \in \{S_1, S_2,\dots,S_{n_e(T)}\}$ равное
числу $\alpha,$ тогда в точке $\S_k \in \sigma(T)$ сушествует
эффект Ефимова.}

\emph{Доказательство.} Допустим, что для $\alpha \in \sigma_d
(H_1)$ существует $S_k \in \{S_1, S_2, \dots, \\ S_{n_e(T)}\}$
такое, что $S_k=\alpha.$ Тогда из представления (12) существенного
спектра существует $\delta > 0$ такое, что $(\alpha-\delta,
\alpha]\cap \sigma_e(T)=\{\alpha\}.$ С другой стороны, имеем
$\omega_n = \alpha+ \beta_n \in \sigma_p(T), \ \forall n \in
\mathbb{N}$ и $\omega_n < \alpha, \forall n \in \mathbb{N}$, где
$\sigma_p(T)-$ множество собственных значений оператора $T.$
Посколько $\lim\limits_{n \rightarrow \infty} \beta_n =0,$ то
$\lim\limits_{n \rightarrow \infty} \omega_n =\S_k.$ Тогда в точке
$\S_k \in \sigma_e(T)$ слева существует эффект Ефимова.
$\blacksquare$

\textbf{Замечание 4.2.} \textit{Пусть $0 < |\sigma_d(H_1)| =p$ и
$|\sigma_d(H_2)|= \infty.$\\
 а)Тогда может быть, что в модели (9)
отсуствует эффект Ефимова. Например, в случае $\sigma_e(H_1)=
\sigma_e(H_2)=[0,1], \ \sigma_d(H_1)=\{-1\}$ и $\sigma_d(H_2) =
\{-\frac 2 n\}_{n \in \mathbb{N}}$ имеем $\sigma_0=[0,2], \
\sigma_1=[-1,0], \ \sigma_2=[-2,0]\cup B,$ где $B \subset [-\frac
2 3, 1].$ Следовательно, $\sigma_e(T)=[-2,2]$ и
$\sigma_d(T)=\{-3\}.$ \\
б) В модели (9) эффект Ефимова появиться не больше чем в $p$
точках.}

\textbf{Теорема 4.6.} \textit{Пусть $|\sigma_d(H_1)| = \infty$ и
$0<|\sigma_d(H_2)| = q.$ Если для $\beta \in \sigma_d(H_2)$
найдется $S_k \in \{S_1, S_2,\dots,S_{n_e(T)}\}$ равное числу
$\beta,$ тогда в точке $\S_k \in \sigma(T)$ сушествует эффект
Ефимова.}

\textbf{Теорема 4.7.} \textit{Пусть $|\sigma_d(H_1)|
=|\sigma_d(H_2)|= \infty.$ Тогда  в каждой точке $S_k \in \{S_1,
S_2,\dots,S_{n_e(T)}\}$ сушествует эффект Ефимова.}

\emph{Доказательство.} Допустим, что $|\sigma_d(H_1)|
=|\sigma_d(H_2)|= \infty$ . Из представления (12) существенного
спектра $\sigma_e(T)$ получаем, что эффект Ефимова может
существовать в точках $S_k, k \in \{1,2,\dots,n_e(T)\}.$ С другой
стороно, для каждого $\S_k$ найдется  номер $\j_k$ такая, что
$\S_k=\alpha_(j_k)$  или $\S_k=\beta_(j_k).$ Пусть имеет  место
равенство $\S_k=\alpha_(j_k).$ Очевидно,  что  $\omega_n =
\alpha_(j_k)+ \beta_n \in \sigma_p(T), \ \forall n \in
\mathbb{N}.$ Имеем $\omega_n < \alpha_(j_k), \forall n \in
\mathbb{N}$ и $\lim\limits_{n \rightarrow \infty} \omega_n
=\alpha_(j_k)=\S_k.$ Значит в точке $\S_k \in \sigma(T)$  слева
существует эффект Ефимова. Анолагично рассматривается случай
$\S_k=\beta_(j_k).$ $\blacksquare$

\textbf{Замечание 4.3.} \textit{Пусть
$|\sigma_d(H_1)|=|\sigma_d(H_2)|= \infty.$ Тогда возможно, что в
модели (9) сушествует эффект Ефимова только в единственной точке
$\omega_0=min\{\alpha_1, \beta_1\},$ т.е. $n_e(T)=1.$ Например, в
случае $\sigma_e(H_1)=\sigma_e(H_2)= [0,1],$
 $\sigma_d(H_1)=\{-\frac 4 n\}_{n \in \mathbb{N}}$ и
 $\sigma_d(H_2)=\{-\frac 5 n\}_{n \in \mathbb{N}}$ имеем
 $\sigma_0=[0,2],$ $\sigma_1=[-4,1]$ и $\sigma_2=[-5,1].$
Следовательно, получаем, что $\sigma_e(T)=[-5,2],\,
\sigma_d(T)\subset [-9, -5)$ и
$\omega_n=-5-\frac 4 n \in \sigma_d(T), n \in \mathbb{N}.$}\\

\textbf{5. Примеры}\\

    Изучим дискретный спектр и эффект Ефимова в модели Фридрихса в
конкретных примерах.

{\bf Случай 1}: $|\sigma_p(U)|= \infty,$ \emph{где} $\sigma_p(U)-$
\emph{множество собственных значений оператора} $U.$ Рассмотрим
последовательность $r_n(x)$ непрерывных функций на $[0,1]$
заданных по формуле
$$
r_n(x) = \left\{ \begin{array}{llll} \frac {3} {h_k}(x-p_k),\quad
\mbox {если} \quad x \in [p_{k}, p_{k}+\frac{h_k}{3}) \\
 1, \quad \mbox{если}  \quad x \in [p_k+\frac{h_k} 3, p_{k+1}-\frac{h_k} 3) \\
-\frac 3 {h_k}(x-p_{k+1}), \quad \mbox{если} \quad x \in
[p_{k+1}-\frac{h_k}3, p_{k+1}] \\
0, \quad \mbox{если} \quad x \notin [p_k, p_{k+1}]\\
\end{array}\right. $$
где $p_{1}=0 ,\ p_{2}= \frac {1} {2} ,\ p_{n+1}=p_{n} + \frac {1}
{2^n}, n \in \mathbb{N},$ $h_{k} =  \frac {1} {2^k}, \ k \in
\mathbb{N}.$

Пусть функция потенциала $u(x)$ в модели (6) задана по формуле
$$u(x)=\sum\limits_{k \in \mathbb{N}}a_k r_k(x), \eqno (13)$$
где $a_n > 0, n \in \mathbb{N}$ и $a_1 = 1, \sum\limits_{n \in
\mathbb{N}} a_n < \infty.$

Ряд (13) равномерно сходится на отрезке $[0,1].$ Следовательно,
каждая функця $u(x)$, заданная по формуле (13), является
неотрицательной и непрерывной на $[0,1],$ при этом $u_{min}=0, $
$u_{max}=1.$

Предположим, что ядро $k(x,s)$ интегрального оператора $K:L_2[0,1]
\rightarrow L_2[0,1]$ задано равенством

$$k(x,s)= \sum\limits_{k \in \mathbb{N}} \lambda_n \varphi_n(x) \varphi_n(s), \eqno (14)$$
где
$$\varphi_n(x)= \frac 3 {\sqrt{h_n}} \ q_n(x), \eqno (15) $$

$$
q_n(x) = \left\{ \begin{array}{lll} \frac 6 {h_n}(x-p_n -
\frac{h_n} 3),\quad
\mbox{если} \quad x \in [p_{n} + \frac{h_n} 3, p_n+\frac{h_n} 2] \\
-\frac 6 {h_n}(x-p_{n+1} +
\frac{h_n} 3), \quad \mbox{если}  \quad x \in (p_n+\frac{h_n} 2, p_{n+1}-\frac{h_n} 3]\\
0, \quad \mbox{если} \quad x \notin [p_n +\frac{h_n} 3, p_{n+1} - \frac{h_n} 3]\\
\end{array}\right. ,
$$

Здесь последовательность $\lambda_n$ удовлетворяет условиям
$\lambda_n >0, n \in \mathbb{N}$ и $\sum\limits_{k \in \mathbb{N}}
2^k \lambda_k < \infty.$ Например, $\lambda_n = \frac 1 {3^n}$
(или $\lambda_n = \frac 1 {a^n}, a > 2$).

Последовательность непрерывных функций $\varphi_n(x)$ (15)
является ортонормированной в $L_2[0,1].$ Ряд (14) равномерно
сходится на квадрате $[0,1] \times [0,1].$ Интегральный оператор
$K$, заданный ядром $k(x,s)$ (14)  является самосопряженным и
положительным в $L_2[0,1].$

Допустим, что для некоторого натурального $n_0$ выполняется
условие $a_n < \lambda_n, \ \forall n \geq n_0.$ Тогда  в модели
Фридрихса ниже существенного спектра существует эффект Ефимова,
так как

$$(H\varphi_n)(x)= (U\varphi_n)(x)-(K \varphi_n)(x)= (a_n -
\lambda_n)\varphi_n(x), \ n \in \mathbb{N},$$ т.е. $\omega_{p_0},
\omega_{p_0+1},\dots, \omega_{p_0+k}, \dots$ являются
отрицательными собственными значениями оператора $U-K,$ где
$\omega_{n}=a_{n}-\lambda_{n}$ и $\lim\limits_{k \rightarrow
\infty} \omega_{p_0+k}=0.$ (Например, в случае $a_n=\frac 1
{5^{n-1}}, \ n \in \mathbb{N}$ и $\lambda_n = \frac 1 {3^n}\ n \in
\mathbb{N}$ числа $\omega_1, \omega_2$ и $\omega_3$ являются
положительными собственными значениями оператора $H=U-K,$ а числа
$\omega_4, \omega_5,\dots, \omega_n,\dots$ являются отрицательными
собственными значениям оператора $U-K,$, где
$\omega_{n}=a_n-\lambda_{n}$ т.е. $\sigma_d
(H-K)=\{\omega_k\}_{k=4}^{\infty}$).

{\bf Случай 2}: $|\sigma_p(U)|=0.$ Пусть $u(x)-$ произвольная
убывающая непрерывная неотрицательная функция на $[0,1]$, которая
удовлетворяет условию $u^{-1}(\{0\})\cap[0,1] \neq \emptyset.$
Положим $$q_n =q_n (u)=\sup \limits_{x \in [p_n,
p_{n+1}]}{u(x)},\,\ n \in \mathbb{N}, $$ где последовательность
$p_k$ определена в примере 1. Тогда $q_n \geq q_{n+1},\,\ n \in
\mathbb{N}$ и $\lim\limits_{n \rightarrow \infty}q_n=0.$ Если
 $\{q_n\}_{n \in \mathbb{N}}\in l_2,$ тогда всегда найдется ядро
 $k(x,s),$ при котором в модели Фридрихса
сушествует эффект Ефимова, так как в пространстве $L_2[0,1]$ мы
можем взять последовательность ортонормированных функций
$$\nu_n(x) = 2^{\frac {{n+1}} {2}} sin \xi_n(x),$$
где
$$
\xi_k(x) = \left\{ \begin{array}{ll} \frac {\pi}
{p_{k+1}-p_k}(x-p_k ),\quad
\mbox{если} \quad x \in [p_k, p_{k+1}] \\
0, \quad \mbox{если} \quad x \notin [p_k , p_{k+1}]\\
\end{array}\right.
$$
и ядро $k(x,s),$ заданное по формуле

$$k(x,s)=\sum\limits_{n \in \mathbb{N}} \lambda_n \nu_n (x) \nu_n(s),$$
где $\lambda_n >0, n \in \mathbb{N},  \sum\limits_{n \in
\mathbb{N}} \lambda_n^2 < \infty$ и для некоторого $n_0 \in
\mathbb{N}$ при этом выполняется неравенство $$q_n<\lambda_n,\,\
\forall n\geq n_0.$$  Тогда в модели (6) существует эффект
Ефимова, т.е. выполняется достаточное условие теоремы 3.4:

$$((U-K)\nu_n, \nu_n) = \int\limits_{p_n}^{p_{n+1}} {u(x)\nu_n^2
(x)dx-\lambda_n}\leq q_n-\lambda_n<0 \,\ \mbox{при всех} \,\ n \geq n_0.$$

{\bf Случай 3}: $|\sigma_d(K)| =\infty$ и $|\sigma_d(U-K)|
<\infty.$ В работах [6-9] доказано, что если функция $u(x)-$
аналитическая в некоторой комплексной окрестности отрезка $[0,1]$
и каждая критическая точка функции $u(x)$ невырождена, то для
произвольного симметрического ядра $k(x,s),$ аналитического  на
$[0,1]^2$ дискретный спектр оператора $U-K$ конечен. Например, в
случае $u(x)=x^2$ для произвольного симметричного ядра $k(x,s),$
аналитического на $[0,1]^2$ в модели Фридрихса отсутствует эффект
Ефимова.

\newpage

\begin{center} \textbf{Сведения об авторе} \end{center}

Ф.И.О: Эшкабилов Юсуп Халбаевич.

Место работы: механико-математический факультет Национального
Университета Узбекистана им М.Улугбека.

Адрес по месту работы: 100174, Узбекистан, г Ташкент, ВУЗгородок.

Домашний адрес: 100139, Узбекистан, г Ташкент, Чиланзар
8371-23-34-42.

Тел: 8371-246-02-30 (раб), 8371-274-67-10 (дом)

e-mail: yusup62@rambler.ru
\end{document}